\documentclass[11pt,a4paper]{amsart}
\usepackage{amsfonts}
\usepackage{amsthm}
\usepackage{amsmath}
\usepackage{amssymb}
\usepackage{amscd}
\usepackage{t1enc}
\usepackage[margin=2.9cm]{geometry}
\usepackage{epstopdf} 
\usepackage{times}
\usepackage{enumerate}
\usepackage{changepage}
\usepackage{mfirstuc}  
\usepackage{tikz}
\usepackage{cite}
\usepackage{booktabs}
\usepackage{tabularx}
\numberwithin{equation}{section}

\theoremstyle{plain}
\newtheorem{Th}{Theorem}[section]

\newcommand{\Q}{\mathbb{Q}}

\newcommand{\Z}{\mathbb{Z}}
\newcommand{\F}{\mathbb{F}}
\newcommand{\C}{\mathbb{C}}

\newcommand{\ord}{\operatorname{ord}}

\begin{document}

\title{On algebraic degrees of inverted Kloosterman sums}
\author{Xin Lin}
\address{Department of Mathematics, Shanghai Maritime University, Shanghai 201306, PR China.}
\email{xlin1126@hotmail.com}

\author{Daqing Wan}
\address{Department of Mathematics, University of California, Irvine, CA
  92697-3875 USA.}
 \email{dwan@math.uci.edu} 
  


\begin{abstract} The study of $n$-dimensional inverted Kloosterman sums was suggested by Katz (1995) who handled 
the case when $n=1$ from complex point of view. For general $n\geq 1$, the $n$-dimensional 
inverted Kloosterman sums were studied from both complex and $p$-adic point of view in our previous paper.  
In this note, we study the algebraic degree of the inverted $n$-dimensional  Kloosterman sum as an algebraic integer.

\end{abstract}

\maketitle

\allowdisplaybreaks[4]
\section{Introduction}\label{sec1}

Let $\F_{q}$ be the finite field of $q$ elements with characteristic $p$ and let $\zeta_p$ be a fixed primitive $p$-th root of unit in $\C$. 
Let $\psi: \F_{q} \rightarrow \mathbb{C}^*$ be the nontrivial additive character defined by 
$$\psi(x) = \zeta_p^{{\rm Tr}(x)},$$
where  ${\rm Tr}$ is the trace map from $\F_q$ to $\F_p$. 
For $b\in\F_{q}^*$ and integer $n\geq 1$, recall that the classical \emph{$n$-dimensional Kloosterman sum} is defined by
\begin{equation*}
	K_n(q,b)=\mathop{\sum_{x_1\cdots x_{n+1}=b, \ x_{i} \in \F^*_{q}}} \psi\left({x_1+\cdots+x_{n+1}}\right)
	=\mathop{\sum_{\ x_{i} \in \F^*_{q}}} \psi\left({x_1+\cdots+x_n+\frac{b}{x_1\cdots x_n}}\right).
\end{equation*}
For each fixed $n$, this is a family of exponential sums, parametrized by the one parameter $b\in \F_q^*$. 
This sum has been studied extensively in the literature. Deligne's deep theorem \cite{Del1980} gives the following 
sharp estimate. 
$$|K_n(q, b)|  \leq (n+1)\sqrt{q}^n.$$
Sperber \cite{Spe1980} determined the $p$-adic slopes for the L-function of the Kloosterman sums. 
These are local results on the number $K_n(q,b)$. 

A global arithmetic problem is to view $K_n(q,b)$ as an algebraic integer in the $p$-th cyclotomic field $\Q(\zeta_p)$ 
and study its degree $\deg K_n(q,b)$ over $\Q$. This degree problem turns out to be rather intricate 
and remains far from well-understood.  A simple result from \cite{Wan1995} says that if ${\rm Tr}(b) \not=0$, then $\deg K_n(q, b)$ is equal to  its maximal possible value $(p-1)/(n+1, p-1)$. 
In particular, this settles the degree problem when $q=p$. If ${\rm Tr}(b)=0$,  then $\deg K_n(q, b)$ can be quite subtle. 
Motivated by \cite{Wan1995}, the degree of a generalized form of $K_n(q, b)$ is presented by Yang\cite{Yan2024} recently.

Katz \cite{Kat1995} raised the question of studying the following 
\emph{inverted $n$-dimensional Kloosterman sum} defined by
\begin{equation*}
	I\!K_n(q, b)=\mathop{\sum_{x_1\cdots x_{n+1}=b, \ x_{i} \in \F^*_{q}}}_{x_1+\cdots+x_{n+1}\neq 0} \psi\left(\frac{1}{x_1+\cdots+x_{n+1}}\right). 
	\end{equation*}
Namely, 
\begin{equation*}	
I\!K_n(q, b) =\mathop{\sum_{\ x_{i} \in \F^*_{q}}}_{x_1+\cdots+\frac{b}{x_1\cdots x_{n}}\neq 0} \psi\left(\frac{1}{x_1+\cdots+x_n+\frac{b}{x_1\cdots x_{n}}}\right).
\end{equation*}
When $n=1$, Katz\cite{Kat1995} obtained a sharp upper bound for $I\!K_1(q, b)$. This result along with paper \cite{Eva1995, Ang1996} motivates the study of Ramanujan graphs.
As we shall see, this inverted sum turns out to be quite interesting in the sense that several  new features occur. 
First, as a complex number, one does not have the expected square root cancellation estimate because there is a non-trivial 
main term. For any positive integer $n$, the sum $I\!K_n(q,b)$ is estimated recently in our previous paper\cite{LW2024}. Two estimates were obtained. 
The first estimate is elementary. The second estimate is much deeper but assumes that $p$ does not divide $n+1$. It is based on 
the theorem of Denef-Loeser \cite{DL1991} for toric exponential sums which in turn depends on Deligne's theorem on the Weil conjectures. 
We state the two estimates below. 

\begin{Th}\label{thm0} 
	Notations as above. We have
	\begin{align*}
		\left|I\!K_n(q, b)+\frac{(q-1)^{n}}{q}\right| \leq q^{\frac{n+1}{2}}.
	\end{align*}
	If $p\nmid n+1$, 
	we further have
	\begin{align*}
		|I\!K_n(q, b)+\frac{(q-1)^{n}-(-1)^{n}(q+1)}{q}|\leq 2n q^{\frac{n}{2}}.
	\end{align*}
	
\end{Th}

In this note, we view the inverted Kloosterman sum $I\!K_n(q, b)$ as an algebraic integer in the $p$-th cyclotomic field 
$\Q (\zeta_p)$ and study its degree as an algebraic integer. This is again quite subtle in general. 
As easily shown later, $\deg I\!K_n(q, b)$ always divides $(p-1)/(n+1, p-1)$. The main result of this note is to prove that the equality holds in the case of prime field $\F_p$.

\begin{Th}\label{thm1} 
Notations as above. 
For all $b \in \F_p^*$, we have 
$$\deg I\!K_n(p, b)  = \frac{p-1}{(n+1, p-1)}.$$ 
\end{Th}

Our method is $p$-adic in nature. We first derive an elementary $p$-adic formula for $I\!K_n(q,b)$ in terms of Gauss sums. 
Then, we apply the classical Stickelberger theorem to identify the $p$-adic main terms that would hopefully distinguish the 
Galois conjugates. This is not always possible for general $q$, but we show it works if $q=p$.  

 For the classical Kloosterman sum $K_n(q, b)$, all the Hodge numbers 
are $1$. However, for the inverted Kloosterman sum $I\!K_n(q,b)$, all the Hodge numbers are $2$, 
except for two of them which are $1$. 
This perhaps explains why the inverted sum $I\!K_n(q,b)$ is more complicated than the classical sum $K_n(q, b)$. 

We conclude this introduction with some further research problems. The first natural problem is to prove degree results 
for the inverted Kloosterman sum $I\!K_n(q,b)$ when $q$ is a proper power of $p$. In this paper, we only treated 
the case $q=p$. It would be interesting to prove some results when $q$ is a proper power of $p$. 
The second problem is to prove degree results for other classical exponential sums which have been 
studied from both complex and $p$-adic point of view, such as the exponential sums in \cite{Wan2004, Wan2021, FW2021, Li2021, CL2022, LC2022, YZ2022, ZF2014}.

\section{A $p$-adic  formula}\label{sec3}

We view the inverted Kloosterman sum as an algebraic integer in the $p$-adic field $\C_p$. 
For this purpose, it suffices to take $\zeta_p$ to be a primitive $p$-th root of unity in $\C_p$. Then, the character 
$\psi: \F_{q} \rightarrow \mathbb{C}_p^*$ becomes a nontrivial additive $p$-adic character. 
Recall that for integer $n\geq 1$ and $b\in\F_{q}^*$, the inverted $n$-dimesnsional Kloosterman sum is defined to be
\begin{equation*}
	I\!K_n(q, b)=\mathop{\sum_{x_1\cdots x_{n+1}=b, \ x_i \in \F_{q}^*}}_{x_1+\cdots+x_{n+1}\neq 0} \psi\left(\frac{1}{x_1+\cdots+x_{n+1}}\right). 
\end{equation*}
The notation $\chi: \F_{q}^* \rightarrow \mathbb{C}_p^*$ denotes a multiplicative $p$-adic character. By the orthogonality of characters, we have
\begin{align}\label{eqs}
	I\!K_n(q, b)=&\frac{1}{q(q-1)}\sum_{\lambda,x_i\in\F^*_q} \sum_{u\in\F_q}\psi\left(u\left(x_1+\cdots+x_{n+1}-\lambda\right)\right) \psi\left(\frac{1}{\lambda}\right)\sum_{\chi}\chi\left(\frac{x_1\cdots x_{n+1}}{b}\right)\nonumber\\
	=&\frac{1}{q(q-1)}\sum_{\lambda\in\F^*_q}\sum_{x_i\in\F^*_q}
	\psi\left(\frac{1}{\lambda}\right)\sum_{\chi}\chi\left(\frac{x_1\cdots x_{n+1}}{b}\right)\nonumber\\
	&+\frac{1}{q(q-1)}\sum_{\lambda\in\F^*_q}\sum_{x_i\in\F^*_q} \sum_{u\in\F_q^*}\psi\left(u\left(x_1+\cdots+x_{n+1}-\lambda\right)\right) 
	\psi\left(\frac{1}{\lambda}\right)\sum_{\chi}\chi\left(\frac{x_1\cdots x_{n+1}}{b}\right)\nonumber\\
	:=&S_1+S_2.
\end{align}
Then
\begin{align}\label{eqs1}
	S_1&=\frac{1}{q(q-1)}\sum_{\lambda\in\F^*_q}\sum_{\chi} \chi^{-1}(b)\psi\left(\frac{1}{\lambda}\right) \sum_{x_i\in\F^*_q}
	\chi(x_1)\cdots \chi(x_{n+1})\nonumber\\
	&=\frac{1}{q(q-1)}\sum_{\lambda\in\F^*_q}\psi\left(\frac{1}{\lambda}\right)
	\sum_{\chi} \chi^{-1}(b)\left(\sum_{x\in\F^*_q}\chi(x)\right)^{n+1}\nonumber\\
	&= -\frac{(q-1)^{n}}{q}. 
\end{align}
Define the Gauss sum $G(\chi)$ by 
$$G(\chi) = \sum_{x \in \F_q^*} \chi(x)\psi(x).$$
Note that if $\chi$ is trivial, then $G(\chi)=-1$. We have 
\begin{align}\label{eqs2}
	S_2&=\frac{1}{q(q-1)}\sum_{\lambda,u\in\F^*_q}\sum_{\chi}\chi^{-1}(b)\sum_{x_i\in\F^*_q}
	\chi(x_1)\psi(ux_1)\cdots\chi(x_{n+1})\psi(ux_{n+1})\psi(-u\lambda)\psi\left(\frac{1}{\lambda}\right)\nonumber\\
	&=\frac{1}{q(q-1)}\sum_{\lambda,u\in\F^*_q}\sum_{\chi}\chi^{-1}(b){\chi^{-(n+1)}}(u)
	\psi(-u\lambda)\psi\left(\frac{1}{\lambda}\right)
	G(\chi)^{n+1}\nonumber\\
	&=\frac{1}{q(q-1)}\sum_{\chi}\chi^{-1}(b)\left(\sum_{\lambda\in\F^*_q}
	{\chi^{-(n+1)}}\left(-\frac{1}{\lambda}\right)\psi\left(\frac{1}{\lambda}\right)\right)
	G(\chi^{-(n+1)}) G(\chi)^{n+1}\nonumber\\
	&=\frac{1}{q(q-1)}\sum_{\chi}\chi^{-1}(b) \chi^{n+1}(-1)
	G(\chi^{-(n+1)})^2 G(\chi)^{n+1}\nonumber\\	
	&=\frac{1}{q(q-1)}\sum_{\chi}\chi^{n+1}(-1) \chi^{-1}(b)
	G(\chi^{-(n+1)})^2 G(\chi)^{n+1}\nonumber\\	
	&=\frac{1}{q(q-1)}\left((-1)^{n+1}+ \sum_{\chi\not=1}\chi^{n+1}(-1) \chi^{-1}(b)
	G(\chi^{-(n+1)})^2 G(\chi)^{n+1} \right).
	\end{align}
It follows that
\begin{align}\label{eqs3}
	I\!K_n(q, b)&=\frac{1}{q(q-1)}\left(-(q-1)^{n+1}+(-1)^{n+1}+ \sum_{\chi\not=1}\chi^{n+1}(-1) \chi^{-1}(b)
	G(\chi^{-(n+1)})^2 G(\chi)^{n+1} \right).	
	\end{align}

\section{Proof of Theorem}

We now turn to the proof of our theorem on the degree of the inverted Kloosterman sum. 
At this point, we do not assume that $q=p$ yet. 

The Galois group of $\Q(\zeta_p)$ over $\Q$ is isomorphic to the cyclic group $\F_p^*$ of order $p-1$. 
For $a\in \F_p^*$, let $\sigma_a$ denote the automorphism such that $\sigma_a (\zeta_p) = \zeta_p^a$. One checks that 
$$\sigma_a(I\!K_n(q, b)) =\mathop{\sum_{x_1\cdots x_{n+1}=b, \ x_i \in \F_{q}^*}}_{x_1+\cdots+x_{n+1}\neq 0} \psi\left(\frac{a}{x_1+\cdots+x_{n+1}}\right)= I\!K_n(q, ba^{-(n+1)}).$$
Let $H= \{ a\in \F_p^* | a^{(n+1, p-1)} =1\}$. This is a cyclic subgroup of order $(n+1, p-1)$ in $\F_p^*$. 
It is clear that $H$ acts trivially on $I\!K_n(q,b)$, namely, $I\!K_n(q,b)$ is in the fixed field $\Q(\zeta_p)^H$. 
By Galois theory, we deduce that $\deg I\!K_n(q, b)$ divides $(p-1)/(n+1, p-1)$, giving the desired upper bound for $\deg I\!K_n(q, b)$. 
Non-trivial lower bound is harder to obtain. One needs to show that 
the set $\{ I\!K_n(q, ba^{-(n+1)}): a\in \F_p^*\}$ contains many different elements. This can be difficult to do in general. 
Our strategy is to find a $p$-adic asymptotic formula for $I\!K_n(q, ba^{-(n+1)})$ and shows that its main 
term contains many different elements as $a$ varies in $\F_p^*$. This will produce a non-trivial lower bound.

For $a\in \F_p^*$, equation (\ref{eqs3}) implies 
\begin{align}\label{eqs4}
	&I\!K_n(q, b) - I\!K_n(q, ba^{-(n+1)}) =\nonumber \\ 
	&\frac{1}{q(q-1)}\left(\sum_{\chi \not=1}\chi^{n+1}(-1)(\chi^{-1}(b) -\chi^{-1}(ba^{-(n+1)}))
	G(\chi^{-(n+1)})^2 G(\chi)^{n+1} \right).	
	\end{align}
We would like to know when this difference of Galois conjugates is non-zero. For this  purpose, we first analyze 
each term when the non-trivial character $\chi$ varies. 
If $\chi^{n+1}=1$, then 
$$\chi^{-1}(ba^{-(n+1)}) = \chi^{-1}(b) \chi^{n+1}(a) = \chi^{-1}(b).$$
This shows that all terms in equation (\ref{eqs4}) with $\chi^{n+1}=1$ are zero, and hence they cannot distinguish the Galois conjugates of $I\!K_n(q, b)$. 
One needs to find other terms with $\chi^{n+1}\not =1$ and with minimal $p$-adic valuation. It is not clear how 
to identify such $p$-adic main terms in general, but the Stickelberger theorem is useful in this regard.  
	
Let $\pi$ be the $(p-1)$th root of $-p$ such that $\pi\equiv\psi(1)-1(\bmod\ \pi^2)$. Then one has $\ord_p(\pi)=1/(p-1)$.
Let $\omega: \F_q^* \rightarrow \C_p^*$ be the Teichm\"uller character. This is a primitive character of order $q-1$. 
Any character $\chi:  \F_q^* \rightarrow \C_p^*$ 
can be written as $\omega^{-m}$ for a unique integer $0\leq m \leq q-2$. The character $\chi^{n+1}$ is non-trivial 
if and only if $(n+1)m$ is not divisible by $q-1$. The Stickelberger theorem says that 
\begin{align}\label{eq1}
	v_p(G(\omega^{-m})) =  \frac{\sigma(m)}{p-1}
	\quad \text{or} \quad
	v_{\pi}(G(\omega^{-m})) =  \sigma(m),
\end{align}
where $\sigma(m)$ is the sum of $p$-digits in the base $p$ expansion of $m$, and $v_p(z)$ denotes the $p$-adic 
valuation of $z$ normalized such that $v_p(p)=1$. 

We now assume that $q=p$ and denote $I\!K_n(p,b)$ as $I\!K_n(b)$. Then, $\chi = \omega^{-m}$ for a unique integer $0\leq m \leq p-2$ and 
the Stickelber theorem simplifies to 
$$v_p(G(\omega^{-m})) =  \frac{\sigma(m)}{p-1} = \frac{m}{p-1}.$$
We only need to consider those $0\leq m\leq p-2$ such that $(n+1)m$ is not 
divisible by $p-1$. In this case, 
equation (\ref{eqs4}) shows that 
\begin{align}\label{eqs5}
	&p(p-1)\left(I\!K_n(b) - I\!K_n(ba^{-(n+1)})\right)=\nonumber \\
	&\left(\mathop{\sum_{1\leq m\leq p-2}}_{(p-1)\nmid (n+1)m}(-1)^{m(n+1)}\omega^m(b)(1- \omega^m(a^{-(n+1)}))
	G(\omega^{m(n+1)})^2 G(\omega^{-m})^{n+1}\right).	
	\end{align}
	
If $n+1$ is divisible by $p-1$, it is clear that $I\!K_n(b)=I\!K_n(ba^{-(n+1)})$ for all $a\in \F_p^*$ and thus 
$\deg I\!K_n(b)=1 =(p-1)/(n+1, p-1)$. The theorem is trivially true. In the following, we assume that $n+1$ is not divisible by $p-1$ and $1\leq m\leq p-2$ with $(p-1)\nmid (n+1)m$. By the Stickelberger theorem, the $p$-adic valuation 
$$V(m): = v_p(G(\omega^{m(n+1)})^2 G(\omega^{-m})^{n+1}) = \frac{(n+1)m + 2 \{-(n+1)m\}_{p-1}}{p-1},$$	
where $\{x\}_{p-1}$ denotes the smallest non-negative residue of $x$ modulo $(p-1)$. We shall find all the terms with lower $p$-adic valuations in different cases.
~\\

\emph{Case I: $p-1>n+1$}

We first consider the case $p-1>n+1$ and 
there is at least one term on the right side of equation (\ref{eqs5}), say $m=1$. 
Since $1\leq m \leq p-2$ and $(n+1)m$ is not divisible by $p-1$, there is a unique integer $0 \leq i\leq n$ such that 
$$(p-1)i < (n+1)m < (p-1)(i+1).$$ 
It follows that
\begin{align*}
	V(m) = \frac{(n+1)m + 2((p-1)(i+1)-(n+1)m)}{p-1} = \frac{2(p-1)(i+1)-(n+1)m}{p-1}
\end{align*} 
and
\begin{align}\label{eq2}
	i+1< V(m) = \frac{(p-1)(i+1) +(p-1)(i+1)-(n+1)m}{p-1} <i+2.
\end{align}
In order for $V(m)$ to be minimal, it is necessary to have  $i=0$. In this case, $m$ satisfies 
$$0< (n+1)m < p-1.$$
Write 
$$p-1= k +(n+1)h, \ 1\leq k \leq n+1.$$ 
Then, $m$ can takes all the values in $\{1, 2, \cdots, h\}$ and one checks that 
$$V(m) = 1 + \frac{k +(n+1)(h-m)}{p-1}.$$
These valuations are distinct as $m$ varies in $\{1, 2, \cdots, h\}$. In fact, 
$$1<V(h) < V(h-1) < \cdots <V(1) < 2.$$
Since $(n+1)m$ lies strictly between $0$ and $p-1$, it follows that $(n+1)m$ is not divisible by $p-1$.  
This implies that for each $1\leq m \leq h$, the character $\omega^{-m(n+1)}$ is non-trivial. 
Reducing equation (\ref{eqs5}), we obtain 
\begin{align}\label{eqs6}
	p(p-1)\left(I\!K_n(b) - I\!K_n(ba^{-(n+1)})\right) \equiv 
	\mathop{\sum_{1\leq m\leq h}} (1- \omega^m(a^{-(n+1)}))
	A_m(b) p^{V(m)} \left(\bmod\ p^{2}\right), 	
	\end{align}
where each $A_m(b)$ is a unit in $\Z_p[\pi]$. 
For $1\leq m \leq h$, it is clear that 
$$1-\omega^m(a^{-(n+1)}) \equiv 1 - a^{-m(n+1)} \left(\bmod\ p\right).$$ 
Let $a\in \F_p^*$ such that $a^{(n+1, p-1)} \not=1$. The coefficient 
$1- a^{-m(n+1)}$ is not zero in $\F_p$ for $m=1$. Let $m^*$ be the largest integer in the interval $[1, h]$ such that 
$1- a^{-m(n+1)}$ is not zero in $\F_p$. Then, 
\begin{align}\label{eqs7}
	&p(p-1)\left(I\!K_n(b) - I\!K_n(ba^{-(n+1)})\right) \nonumber \\ 
	&\equiv (1- \omega^{m^*}(a^{-(n+1)}))
	A_m(b) p^{V(m^*)} \left(\bmod\ p^{V(m^*)+ \frac{1}{p-1}}\right) \nonumber \\ 	
	&\not\equiv 0  \left(\bmod\ p^{V(m^*)+ \frac{1}{p-1}}\right). \nonumber
	\end{align}
This shows that $I\!K_n(b) - I\!K_n(ba^{-(n+1)}) \not=0$. In summary, for $a \in \F_p^*$, we have proved that $I\!K_n(b) - I\!K_n(ba^{-(n+1)}) =0$ if and only if $a^{(n+1, p-1)} =1$. It follows that 
$\deg I\!K_n(b) = (p-1)/(n+1, p-1)$. 
~\\

\emph{Case II: $\frac{1}{2}(n+1)<p-1<n+1$}

In this case, there is a unique integer $1 \leq i\leq n-1$ such that 
\begin{align*}
	(p-1)i < (n+1)m < (p-1)(i+1).
\end{align*}
Since we assume $p-1<n+1<2(p-1)$, the integer $m$ is uniquely determined by $i$ and $i$ is also bounded by $m$.
As shown in equation (\ref{eq2}), 
the minimal $V(m)$ corresponds to $i=1$ and one has 
$$(p-1)<(n+1)m<2(p-1),$$ 
which implies $m=1$ and
\begin{align*}
	2<V(1)=2+\frac{2(p-1)-(n+1)}{p-1}<3.
\end{align*}
For $m\geq 2$, the restriction $(p-1)i < (n+1)m < (p-1)(i+1)$ implies that $i\geq 2$ and thus $V(m)>3$.
Combining with equation (\ref{eqs5}), one obtains
\begin{align*}
	p(p-1)\left(I\!K_n(b) - I\!K_n\left(ba^{-(n+1)}\right)\right)
	\equiv
	\left(1- \omega\left(a^{-(n+1)}\right)\right)p^{V(1)}C_1(b) \left(\bmod\ p^3\right),
\end{align*}
where $C_1(b)$ is a unit in $\Z_p[\pi]$. 
Then for $a \in \F_p^*$, we have proved that
$I\!K_n(b) = I\!K_n(ba^{-(n+1)})$ if and only if $a^{(n+1, p-1)} =1$.
~\\

\emph{Case III: $p-1<\frac{1}{2}(n+1)$}

By the Stickelberger theorem, the $p$-adic valuation 
\begin{equation}\label{eq3}
	W(m): = v_{\pi}(G(\omega^{m(n+1)})^2 G(\omega^{-m})^{n+1}) = (n+1)m + 2 \{-(n+1)m\}_{p-1}.
\end{equation}
We shall determine when $W(m)$ reaches its minimum as $1\leq m\leq p-2$. Denote $h=\{n+1\}_{p-1}$. Since we assumed that $n+1$ is not divisible by $p-1$, one has $1\leq h\leq p-2$. For $m=1$, it can be deduced from equation $(\ref{eq3})$ that $W(1)=n+1+2(p-1-h)$. For $2\leq m\leq p-2$, there is a unique integer $i\geq 4$ such that 
\begin{align*}
	(p-1)i < (n+1)m < (p-1)(i+1).
\end{align*}
Then one has
\begin{align}\label{eq4}
	W(m)-W(1)=(m-1)(n+1)-2(p-1)+2\left((i+1)(p-1)-m(n+1)\right)+2h.
\end{align}
 Since we assume $n+1>2(p-1)$ in this case, equation $(\ref{eq4})$ implies that $W(m)>W(1)+4$ for $2\leq m\leq p-2$. Combining with equation (\ref{eqs5}), one obtains
 \begin{align*}
 	p(p-1)\left(I\!K_n(b) - I\!K_n\left(ba^{-(n+1)}\right)\right)
 	\equiv
 	\left(1- \omega\left(a^{-(n+1)}\right)\right)\pi^{W(1)}C_2(b) \left(\bmod\ \pi^{W(1)+1}\right),
 \end{align*}
 where $C_2(b)$ is a unit in $\Z_p[\pi]$. 
 Then for $a \in \F_p^*$, we have proved that
 $I\!K_n(b) = I\!K_n(ba^{-(n+1)})$ if and only if $a^{(n+1, p-1)} =1$.
The theorem is proved.


\nocite{*}
\bibliographystyle{amsalpha}
\bibliography{ref}


\end{document}